
\newcount\secno
\newcount\prmno
\newif\ifnotfound
\newif\iffound

\def\namedef#1{\expandafter\def\csname #1\endcsname}
\def\nameuse#1{\csname #1\endcsname}

\long\def\ifundefined#1#2#3{\expandafter\ifx\csname
  #1\endcsname\relax#2\else#3\fi}
\def\hwrite#1#2{{\let\the=0\edef\next{\write#1{#2}}\next}}

\toksdef\ta=0 \toksdef\tb=2
\long\def\leftappenditem#1\to#2{\ta={\\{#1}}\tb=\expandafter{#2}%
                                \edef#2{\the\ta\the\tb}}
\long\def\rightappenditem#1\to#2{\ta={\\{#1}}\tb=\expandafter{#2}%
                                \edef#2{\the\tb\the\ta}}

\def\lop#1\to#2{\expandafter\lopoff#1\lopoff#1#2}
\long\def\lopoff\\#1#2\lopoff#3#4{\def#4{#1}\def#3{#2}}

\def\ismember#1\of#2{\foundfalse{\let\given=#1%
    \def\\##1{\def\next{##1}%
    \ifx\next\given{\global\foundtrue}\fi}#2}}

\def\section#1{\vskip1truecm
               \global\def\currenvir{section}
               \global\advance\secno by1\global\prmno=0
               {\bf \number\secno. {#1}}
               \smallskip}

\def\subsection{\global\def\currenvir{subsection}
                \global\advance\prmno by1
                \smallskip \ind{ (\number\secno.\number\prmno) }}
\def\subsec#1{\global\def\currenvir{subsection}
                \global\advance\prmno by1
                \smallskip{(\number\secno.\number\prmno)\ {\it #1}\ind}}

\def\proclaim#1{\global\advance\prmno by 1
                {\bf #1 \the\secno.\the\prmno$.-$ }}

\long\def\th#1 \enonce#2\endth{%
   \medbreak\proclaim{#1}{\it #2}\global\def\currenvir{th}\smallskip}

\def\bib#1{\rm #1}
\long\def\thr#1\bib#2\enonce#3\endth{%
\medbreak{\global\advance\prmno by 1\bf#1\the\secno.\the\prmno\ 
\bib{#2}$\!.-$ } {\it
#3}\global\def\currenvir{th}\smallskip}
\def\rem#1{\global\advance\prmno by 1
{\it #1} \the\secno.\the\prmno$.-$ }


\def\isinlabellist#1\of#2{\notfoundtrue%
   {\def\given{#1}%
    \def\\##1{\def\next{##1}%
    \lop\next\to\za\lop\next\to\zb%
    \ifx\za\given{\zb\global\notfoundfalse}\fi}#2}%
    \ifnotfound{\immediate\write16%
                 {Warning - [Page \the\pageno] {#1} No reference found}}%
                \fi}%
\def\ref#1{\ifx\labellist\empty{\immediate\write16
                 {Warning - No references found at all.}}
               \else{\isinlabellist{#1}\of\labellist}\fi}

\def\newlabel#1#2{\rightappenditem{\\{#1}\\{#2}}\to\labellist}
\def\labellist{}

\def\label#1{%
  \def\given{th}%
  \ifx\given\currenvir%
    {\hwrite\lbl{\string\newlabel{#1}{\number\secno.\number\prmno}}}\fi%
  \def\given{section}%
  \ifx\given\currenvir%
    {\hwrite\lbl{\string\newlabel{#1}{\number\secno}}}\fi%
  \def\given{subsection}%
  \ifx\given\currenvir%
    {\hwrite\lbl{\string\newlabel{#1}{\number\secno.\number\prmno}}}\fi%
  \def\given{subsubsection}%
  \ifx\given\currenvir%
  {\hwrite\lbl{\string%
    \newlabel{#1}{\number\secno.\number\subsecno.\number\subsubsecno}}}\fi
  \ignorespaces}

\newwrite\lbl

\def\openall{\openout\lbl=\jobname.lbl}

\newread\testfile
\def\lookatfile#1{\openin\testfile=\jobname.#1
    \ifeof\testfile{\immediate\openout\nameuse{#1}\jobname.#1
                    \write\nameuse{#1}{}
                    \immediate\closeout\nameuse{#1}}\fi%
    \immediate\closein\testfile}%

\newlabel{not}{1.1}
\newlabel{w2}{1.3}
\newlabel{sw}{1.4}
\newlabel{main}{1.6}
\newlabel{tl}{2.3}
\newlabel{nsc}{2.4}
\newlabel{ver}{2.5}
\newlabel{gen}{2.7}
\newlabel{sdc}{4.3}
\newlabel{span}{4.4}
\newlabel{conj}{4.5}
\newlabel{span2}{4.6}
\newlabel{cases}{4.7}


\magnification 1250
\pretolerance=500 \tolerance=1000  \brokenpenalty=5000
\mathcode`A="7041 \mathcode`B="7042 \mathcode`C="7043
\mathcode`D="7044 \mathcode`E="7045 \mathcode`F="7046
\mathcode`G="7047 \mathcode`H="7048 \mathcode`I="7049
\mathcode`J="704A \mathcode`K="704B \mathcode`L="704C
\mathcode`M="704D \mathcode`N="704E \mathcode`O="704F
\mathcode`P="7050 \mathcode`Q="7051 \mathcode`R="7052
\mathcode`S="7053 \mathcode`T="7054 \mathcode`U="7055
\mathcode`V="7056 \mathcode`W="7057 \mathcode`X="7058
\mathcode`Y="7059 \mathcode`Z="705A
\def\spacedmath#1{\def\packedmath##1${\bgroup\mathsurround =0pt##1\egroup$}
\mathsurround#1
\everymath={\packedmath}\everydisplay={\mathsurround=0pt}}
\def\nospacedmath{\mathsurround=0pt
\everymath={}\everydisplay={} } \spacedmath{2pt}
\def\qfl#1{\buildrel {#1}\over {\longrightarrow}}
\def\phfl#1#2{\normalbaselines{\baselineskip=0pt
\lineskip=10truept\lineskiplimit=1truept}\nospacedmath\smash {\mathop{\hbox to
8truemm{\rightarrowfill}}
\limits^{\scriptstyle#1}_{\scriptstyle#2}}}
\def\hfl#1#2{\normalbaselines{\baselineskip=0truept
\lineskip=10truept\lineskiplimit=1truept}\nospacedmath\smash{\mathop{\hbox to
12truemm{\rightarrowfill}}\limits^{\scriptstyle#1}_{\scriptstyle#2}}}

\def\iso{\vbox{\hbox to .8cm{\hfill{$\scriptstyle\sim$}\hfill}
\nointerlineskip\hbox to .8cm{{\hfill$\longrightarrow $\hfill}} }}

\def\sdir_#1^#2{\mathrel{\mathop{\kern0pt\oplus}\limits_{#1}^{#2}}}
\def\pprod_#1^#2{\raise
2pt \hbox{$\mathrel{\scriptstyle\mathop{\kern0pt\prod}\limits_{#1}^{#2}}$}}
\font\tengrit=cmmib10
\font\sevengrit=cmmib7
\font\fivegrit=cmmib5
\newfam\gritfam
\textfont\gritfam=\tengrit \scriptfont\gritfam=\sevengrit
  \scriptscriptfont\gritfam=\fivegrit
  \def\grit{\fam\gritfam\tengrit}%
\mathchardef\ll="7860

\font\eightrm=cmr8         \font\eighti=cmmi8
\font\eightsy=cmsy8        \font\eightbf=cmbx8
\font\eighttt=cmtt8        \font\eightit=cmti8
\font\eightsl=cmsl8        \font\sixrm=cmr6
\font\sixi=cmmi6           \font\sixsy=cmsy6
\font\sixbf=cmbx6\catcode`\@=11
\def\eightpoint{%
  \textfont0=\eightrm \scriptfont0=\sixrm \scriptscriptfont0=\fiverm
  \def\rm{\fam\z@\eightrm}%
  \textfont1=\eighti  \scriptfont1=\sixi  \scriptscriptfont1=\fivei
  \def\oldstyle{\fam\@ne\eighti}\let\old=\oldstyle
  \textfont2=\eightsy \scriptfont2=\sixsy \scriptscriptfont2=\fivesy
  \textfont\itfam=\eightit
  \def\it{\fam\itfam\eightit}%
  \textfont\slfam=\eightsl
  \def\sl{\fam\slfam\eightsl}%
  \textfont\bffam=\eightbf \scriptfont\bffam=\sixbf
  \scriptscriptfont\bffam=\fivebf
  \def\bf{\fam\bffam\eightbf}%
  \textfont\ttfam=\eighttt
  \def\tt{\fam\ttfam\eighttt}%
  \abovedisplayskip=9pt plus 3pt minus 9pt
  \belowdisplayskip=\abovedisplayskip
  \abovedisplayshortskip=0pt plus 3pt
  \belowdisplayshortskip=3pt plus 3pt 
  \smallskipamount=2pt plus 1pt minus 1pt
  \medskipamount=4pt plus 2pt minus 1pt
  \bigskipamount=9pt plus 3pt minus 3pt
  \normalbaselineskip=9pt
  \setbox\strutbox=\hbox{\vrule height7pt depth2pt width0pt}%
  \normalbaselines\rm}\catcode`\@=12

\newcount\noteno
\noteno=0
\def\up#1{\raise 1ex\hbox{\sevenrm#1}}
\def\no{n\up{o}\kern 2pt}
\def\note#1{\global\advance\noteno by1
\footnote{\parindent0.4cm\up{\number\noteno}\
}{\vtop{\eightpoint\baselineskip12pt\hsize15.5truecm\noindent
#1}}\parindent 0cm}
\font\san=cmssdc10
\def\ext{\hbox{\san \char3}}
\def\sym{\hbox{\san \char83}}

\def\pc#1{\tenrm#1\sevenrm}
\def\tx{\kern-1.5pt -}
\def\cqfd{\kern 2truemm\unskip\penalty 500\vrule height 4pt depth 0pt width
4pt\medbreak} \def\carre{\vrule height 4pt depth 0pt width 4pt}
\def\virg{\raise
.4ex\hbox{,}}
\def\ind{\par\hskip 1truecm\relax}
\def\indp{\par\hskip 0.5truecm\relax}
\def\moins{\mathrel{\hbox{\vrule height 3pt depth -2pt width 6pt}}}
\def\rond{\kern 1pt{\scriptstyle\circ}\kern 1pt}
\def\Pic{\mathop{\rm Pic}\nolimits}
\def\dim{\mathop{\rm dim}\nolimits}
\def\tl{T_\ell ^{\rm reg}}
\def\ms#1{{\cal M}_{SO_{#1}}}
\def\ls#1{{\cal L}_{SO_{#1}}}
\def\mso{{\cal M}_{SO_r}}
\def\msl{{\cal M}_{SL_r}}
\def\lso{{\cal L}_{SO_r}}
\def\lsl{{\cal L}_{SL_r}}
\def\msp#1{{\cal M}_{{\rm Sp}_{2#1}}}

\def\sp{\mathop{\rm Sp}\nolimits}
\def\ot#1{^{{\scriptscriptstyle\otimes}#1}}
\frenchspacing
\input amssym.def
\input amssym
\vsize = 25truecm
\hsize = 16truecm
\voffset = -.5truecm
\parindent=0cm
\baselineskip15pt
\def\bu{\raise1pt\hbox{$\scriptscriptstyle\bullet$}}
\centerline{\bf  Orthogonal bundles on curves and theta functions}
\smallskip
\smallskip \centerline{Arnaud {\pc BEAUVILLE}} 
\vskip1.2cm

{\bf Introduction}
\ind Let  $C$ be a  curve of genus $g\ge 2$, $G$ an almost  simple
complex Lie group, and ${\cal M}_G$ the moduli space of semi-stable $G$\tx
bundles on $C$. For each component ${\cal M}_G^{\bu}$ of ${\cal M}_G$, the
Picard group is infinite cyclic; its positive generator ${\cal L}_G^{\bu}$ can be
described explicitely as a determinant bundle. 
Then a natural  question, which we will address in this paper
for the classical groups, is to describe the space of ``generalized theta
functions" $H^0({\cal M}_G^{\bu},{\cal L}_G^{\bu})$ and the associated rational 
map
$\varphi_G^{\bu}:{\cal M}_G^{\bu}\dasharrow |{\cal L}_G^{\bu}|_{}^*$. 
\ind The model we have in mind is the case $G=SL_r$. Let $J^{g-1}$ be the
component of the Picard variety of $C$ parameterizing line bundles of degree 
$g-1$; it is isomorphic to the Jacobian of $C$, and carries a canonical theta divisor
$\Theta $ consisting of  line bundles $L$ in $J^{g-1}$ with $H^0(C,L)\not=0$.  
For a general $E\in \msl$, the locus
$$\Theta_E= \{L\in J^{g-1}\ |\ H^0(C,E\otimes L)\not= 0\}$$is in a natural
way a divisor, which belongs to the linear system $|r\Theta|$ on $J^{g-1}$.
We thus obtain a rational map $\vartheta:\msl \dasharrow |r\Theta|$. The main
result of [BNR] is that {\it there exists an isomorphism $|\lsl|^*\iso 
|r\Theta|$ which identifies the rational maps $\varphi_{SL_r}$ and} $\vartheta$.
 This gives a reasonably concrete description of  $\varphi_{SL_r}$, which
allows to get some  information on the behaviour of this map, at least for small
values of $r$ or $g$ (see [B3] for a survey of recent results).
\ind Let us consider now the case $G=SO_r$ with $r\ge 3$. The moduli space 
$\mso$ parametrizes vector bundles $E$ on $C$ of rank $r$ and trivial
determinant,
 with a non-degenerate quadratic form
$q:\sym^2E\rightarrow {\cal O}_C$; it has two components $\mso^+$ and
$\mso^-$. Let $\theta :\mso\dasharrow |r\Theta |$ be the map $(E,q)\mapsto
\Theta _E$. We will see that $\theta $ maps $\mso^+$ and $\mso^-$ into the
subspaces
$|r\Theta|^+ $ and $|r\Theta |^-$  corresponding to
even and odd theta functions respectively. Our main result is:
\smallskip 
{\bf Theorem}$.-$ {\it There are canonical isomorphisms
$|{\cal L}^{\pm}_{SO_r}|^*\iso |r\Theta |^{\pm}$  which identify 
$\varphi _{SO_r}^{\pm}:\mso^{\pm}\dasharrow |{\cal L}^{\pm}_{SO_r}|^*$ with 
the map
$\theta ^{\pm}:\mso^{\pm}\dasharrow |r\Theta |^{\pm}$ induced by}
$\theta $.\smallskip 
\ind This is easily seen to be  equivalent to the fact
that the pull-back map
$\theta ^*:H^0(J^{g-1}, {\cal O}(r\Theta ))^*\rightarrow
H^0(\mso,\lso)$ is an {\it isomorphism}. We will prove that it is injective by
restricting to a small subvariety of
$\mso$ (\S 1). Then we will use the Verlinde formula (\S 2 and 3) to show that the
dimensions are the same. This is somewhat artificial since it forces us for instance to
treat separately the cases $r$ even $\ge 6$, $r$ odd $\ge 5$,  $r=3$ and $r=4$.
 It would be interesting
to find  a more direct proof, perhaps in the spirit of [BNR].
\ind In the last section we consider the same question for the symplectic group. 
Here the theta map does not involve the Jacobian of $C$
but the moduli space 
${\cal N}$ of semi-stable rank 2 vector bundles on $C$ with determinant $K_C$. Let
${\cal L}$ be the determinant bundle on ${\cal N}$. For  $(E,\varphi )$ general in
$\msp{r}$, the reduced subvariety
$$\Delta _E=\{F\in {\cal N}\ |\ H^0(E\otimes F)\not=0\}$$is a divisor on ${\cal N}$,
which belongs to the linear system $|{\cal L}^r|$; this defines a map
$\msp{r}\dasharrow |{\cal L}^r|$
which should coincide, up to a canonical isomorphism, with $\varphi _{\sp_{2r}}$.
This is a particular case of the {\it strange duality} conjecture for the
symplectic group, which we discuss in \S 4 . Unfortunately even this particular case
is not known, except in a few cases that we explain below.
\font\cag=cmbsy10
\section{The moduli space $\hbox{\cag\char77}_{{\bf 
SO_{\grit r}}}$}
\subsection Throughout the paper we fix a complex curve $C$ of genus $g\ge 2$. For
$G$ a  semi-simple complex Lie group, we denote by ${\cal M}_G$ the moduli space of 
semi-stable $G$\tx bundles on $C$. It is a normal projective variety, of dimension
$(g-1)\dim G$. Its connected components are in one-to-one correspondence with the
elements of the group $\pi _1(G)$.\label{not}
\subsection Let us consider the case $G=SO_r\,\ (r\ge 3)$. The  space $\mso$ is
the mo\-duli space of (semi-stable) orthogonal bundles, that is pairs $(E,q)$ where $E$
is a semi-stable\note{By [R], 4.2, an orthogonal bundle $(E,q)$ is semi-stable
if and only if the vector bundle $E$ is semi-stable.} vector bundle of rank
$r$ and trivial determinant and
$q:\sym^2E\rightarrow {\cal O}_C$ a non-de\-ge\-nerate quadratic form.
 The two components
$\mso^+$ and $\mso^-$ are distingui\-shed by the parity of the second
Stiefel-Whitney class $w_2(E,q)\in H^2(C,{\bf Z}/2)\cong {\bf Z}/2$. This class has
the following property (see e.g. [Se], Thm. 2): for every theta-characteristic
$\kappa$ on
$C$ and orthogonal bundle $(E,q)\in\mso$,\global\advance\prmno by1
$$w_2(E,q)\equiv h^0(C,E\otimes \kappa)+rh^0(C,\kappa)\ \hbox{ (mod. 2)
.}\eqno{(\number\secno.\number\prmno)}$$
 \ind  The involution $\iota :L\mapsto K_C\otimes L^{-1} $
of $J^{g-1}$ preserves \label{w2}
$\Theta $, hence lifts to an involution of ${\cal O}_{J^{g-1}}(\Theta )$. We denote 
by $|r\Theta |^+$ and $|r\Theta |^-$ the two corresponding eigenspaces in
 $|r\Theta |$, and by $\theta:\mso\dasharrow |r\Theta | $ the map 
$ (E,q)\mapsto
\Theta _E$.
\th Lemma
\enonce The rational map $\theta :\mso\dasharrow |r\Theta |$ maps $\mso^+$ in
$|r\Theta |^+$ and $\mso^-$ in $|r\Theta |^-$.
\endth\label{sw}
{\it Proof} : For any $E\in\msl$ we have $\iota ^*\Theta _E=\Theta _{E^*}$, so 
$\theta (\mso)$ is contained in the fixed locus $|r\Theta |^+\cup |r\Theta |^-$ of
$\iota ^*$. Since $\mso^{\pm}$ is connected, it suffices to find one element
$(E,q)$ of $\mso^+$ (resp. $\mso^-$) such that $\Theta _E$ is a divisor in $|r\Theta
|^+$ (resp. $|r\Theta |^-$).
 \ind Let $\kappa\in J^{g-1}$ be an even
theta-characteristic of $C$; a symmetric divisor  
$D\in |r\Theta |$ is in $|r\Theta |^+$ (resp. $|r\Theta |^-$) if and only if
 ${\rm mult}_\kappa(D)$ is even (resp. odd) -- see [M1],  \S 2.
 Let $J[2]$ be the 2-torsion subgroup of $\Pic(C)$; we take $E=\alpha
_1\oplus\ldots\oplus\alpha _r$, where $\alpha_1,\ldots ,\alpha_r\in J[2]$ and
$\sum \alpha_i=0$. We endow $E$ with the diagonal quadratic form $q$ deduced
from the isomorphisms $\alpha_i^2\cong {\cal O}_C$. Then $\Theta _E=\Theta
_{\alpha_1}+\ldots +\Theta _{\alpha_r}$. By the Riemann singularity
theorem  the multiplicity  at
$\kappa$ of $\Theta _\alpha$ is $h^0(\alpha\otimes \kappa)$. Thus by (\ref{w2})
$${\rm mult}_\kappa(\Theta _E)=\sum_i h^0(\alpha_i\otimes \kappa)=h^0(E\otimes
\kappa)\equiv w_2(E,q)\ \hbox{ mod. 2}\ .\kern
2truemm\carre$$\vskip-5pt
\subsection Let $\lso$ be the determinant bundle on $\mso$, that is, the pull back
of $\lsl$ by the map $(E,q)\mapsto E$, and let
$\lso^{+}$ and $\lso^-$ be its restrictions to
$\mso^{+}$ and $\mso^-$. It follows from [BLS] that  for $r\not=4$,
$\lso^{\pm}$ generates
$\Pic(\mso^{\pm})$. 
\th Proposition
\enonce The map $\theta ^*:H^0(J^{g-1},{\cal O}(r\Theta
))^*\longrightarrow  H^0(\mso,\lso)$ induced by $\theta :\mso \dasharrow
|r\Theta |$ is an 
isomorphism.
\endth\label{main}
\ind By Lemma (\ref{sw}) $\theta^*$ splits as a  direct sum $(\theta
^+)^*\oplus (\theta ^-)^*$, where
$$(\theta
^{\pm})^*:\bigl(H^0(J^{g-1},{\cal
O}(r\Theta))^{\pm}\bigr)^* \ \hfl{}{}\ H^0(\mso^{\pm},\lso^{\pm})\
.$$
The Proposition
implies that $(\theta ^+)^*$ and $(\theta ^-)^*$ are isomorphisms, and this is
equivalent to the Theorem stated in the introduction.
\smallskip 
{\it Proof of the Proposition} : We will show in \S  3 that the Verlinde formula
gives $$\dim H^0(\mso,\lso)=\dim H^0(J^{g-1},{\cal O}(r\Theta))=r^g\ .$$ It is
therefore sufficient to prove that  $\theta^*$ is injective, or equivalently  that
$\theta (\mso)$ spans the projective space $|r\Theta |$.  We consider again the
orthogonal bundles $(E,q)=\alpha _1\oplus\ldots \oplus\alpha _r$ for $\alpha
_1,\ldots, \alpha _r$ in $J[2]$, $\sum \alpha _i=0$. This bundle has a theta
divisor $\Theta _E=\Theta _{\alpha _1}+\ldots +\Theta _{\alpha _r}$. We claim
 that divisors of this form span $|r\Theta |$. To prove it, let us identify
$J^{g-1}$ with the Jacobian $J$ of $C$ (by choosing a divisor class of
degree
$g-1$). For
$a\in J$, the divisor
$\Theta_a$ is the only element of the linear system $|{\cal O}_J(\Theta )\otimes
\varphi (a)|$, where $\varphi :J\rightarrow \hat J$ is the isomorphism associated 
to the principal polarization of $J$.
Therefore our  assertion  follows from the
following easy lemma:
\th Lemma
\enonce Let $A$ be an abelian variety, $L$ an ample line bundle on $A$, $\hat 
A[2]$ the $2$\tx torsion subgroup of $\Pic(A)$. The multiplication map
$$\sum_{\alpha _1,\ldots ,\alpha _r\in\hat A[2]\atop \alpha _1+\ldots + \alpha
_r=0} H^0(A,L\otimes
\alpha_1)\otimes 
\ldots
\otimes H^0(A,L\otimes \alpha_r)\longrightarrow H^0(A,L^r)$$is surjective.
\endth
{\it Proof} : Let $2_A$ be the multiplication by $2$ in $A$. We have canonical
isomorphisms
$$H^0(A,2_A^*L)\cong \bigoplus_{\alpha \in \hat A[2]} H^0(L\otimes \alpha )
\quad ,\quad H^0(A,2_A^*L^r)\cong \bigoplus_{\beta  \in \hat A[2]}
H^0(L^r\otimes
\beta  )\ ;$$through these isomorphisms the product map
$m_r:H^0(A,2_A^*L)\ot{r}\longrightarrow H^0(A,2_A^*L^r)$ is  the direct sum 
over $\beta \in\hat A[2]$ of the maps 
$$m_r^\beta :\sum_{\alpha _1,\ldots ,\alpha _r\in\hat A[2]\atop \alpha
_1+\ldots +\alpha _r=\beta } H^0(A,L\otimes
\alpha_1)\otimes
\ldots
\otimes H^0(A,L\otimes \alpha_r)\longrightarrow H^0(A,L^r\otimes \beta) \  .$$
Since the line bundle $2_A^*L$ is algebraically equivalent to $L^4$,  the map
$m_r$ is surjective [M2], hence so is $m_r^\beta $ for every $\beta $. The case
$\beta =0$ gives the lemma.\cqfd

\section{The Verlinde formula}
\subsection We keep the notation of (\ref{not}); we denote by $q$ the number of
simple factors of the Lie algebra of $G$  (we are mainly interested
in the case $q=1$).
\ind  To each representation $\rho
:G\rightarrow SL_r$ is attached a line bundle ${\cal L}_\rho $ on ${\cal M}_G$, the
pull back of the determinant bundle on $\msl $ by the morphism ${\cal
M}_G\rightarrow \msl$ associated to $\rho $. The Verlinde
formula expresses the dimension of $H^0({\cal M}_G,{\cal L}_\rho ^k)$, for
each integer $k$, in the form 
$$\qquad \dim H^0({\cal M}_G,{\cal L}_\rho ^k)=N_{k{\grit d_\rho} }(G)\ ,\quad
\hbox{ where}$$
\ind $\bullet$
 ${\grit d_\rho} \in {\bf N}^q$ is the {\it Dynkin index} of $\rho $. For
$q=1$ the number $d_\rho $ is defined and computed in [D], \S 2.  In the
general case the universal cover of $G$ is a product $G_1\times \ldots
\times G_q$ of almost simple factors, and we put
${\grit d_\rho}=(d_{\rho _1},\ldots ,d_{\rho _q})$, where $\rho _i$ is the pull
back of $\rho$ to $G_i$.  
\ind We will need only to know that the Dynkin index is 2 for the standard
representation  of
$SO_r\ (r\ge 5)$, $4$ for that of $SO_3$, and $(2,2)$ for that of $SO_4$.
\ind $\bullet$
 $N_{\ll}(G)$ is an integer depending
on   $G$,  the genus $g$ of $C$, and $ \ll\in {\bf N}^q $. We will now explain how
this number is computed. Our basic reference is [AMW].\par
\subsec{ The simply connected case}
 Let us first consider the case where $G$ is
simply connected and almost simple (that is, $q=1$).
 Let $T$ be a maximal torus of $G$, and $R=R(G,T)$
the corresponding root system (we view the roots of
$G$ as characters of $T$). We denote by  $T_\ell$  the (finite) subgroup of
elements $t\in T$ such that ${\alpha}(t)=1$
  for each {\it long} root $\alpha$, and by 
$T^{\rm{reg}}_\ell $ the subset of  regular  elements  $t\in T_\ell $  (that is,
such that $\alpha(t)\not=1$ for each root $\alpha$). It is stable under the
action of the Weyl group $W$. For $t\in T$, we put
$\displaystyle \Delta (t)=\prod_{\alpha \in R}(\alpha (t)-1)$. Then the
Verlinde formula is $$N_\ell (G)= \sum_{t\in T^{\rm reg}_\ell
/W}\Bigl({|T_\ell |\over \Delta (t)}\Bigr)^{g-1}\ .$$
\subsection\label{tl} This number can be explicitely computed in the following way. 
Let ${\goth t}$ be the Lie algebra of
$T$. The character group $P(R)$ of
$T$ embeds naturally into ${\goth t}^*$.
 We endow ${\goth t}^*$ with the $W$\tx invariant bilinear form $(\ |\ )$ such
that
$(\alpha \,|\,\alpha )=2$ for each long root $\alpha $, and we use this product
to identify
${\goth t}^*$ with ${\goth t}$. Let $\theta $ be the
highest root of $R$; we denote by 
$P_\ell $  the set of dominant weights $\lambda \in P(R)$ such that $(\lambda
\,|\,\theta  )\le \ell $. Let $\rho \in P(R)$ be the half-sum of the positive roots. The
number
$h:=(\rho\,|\,\theta  )+1$ is the {\it dual Coxeter number} of $R$. We have
$|T_\ell |=(\ell +h)^sf\nu $, where $s$ is the rank of $R$,
$f$ the order of the center of $G$, and $\nu $  a number depending on $R$; it is 
equal to $1$ for $R$ of type $D_s$ and to $2$ for $B_s$ ([B2], 9.9).
\ind For
$\lambda \in P_\ell $ we put $\displaystyle t_\lambda =\exp 2\pi
i\ {\lambda+\rho\over \ell +h}\ \cdot $ The map $\lambda \mapsto t_\lambda $ is a
bijection of $P_\ell $ onto $\tl/W $ ([B2], 9.3.c)). For $\lambda \in P_\ell $,
we have $\displaystyle
\alpha (t_\lambda )=\exp 2\pi  i\, {(\alpha\, |\,\lambda +\rho  )\over \ell
+h}$, and therefore
$$\Delta (t_\lambda )=\prod_{\alpha \in R_+}4\,\sin^2 \pi\; {(\alpha\,
|\,\lambda +\rho )\over \ell +h}\ \cdot $$
\subsec{The non-simply connected  case}\label{nsc}We now give the formula for a
general almost simple group, following [AMW]. 
\ind Let $Z$ be the center of $G$. An element $t$ of $T$ belongs to $Z$ if and
only if $\alpha (t)=1$ for all $\alpha \in R$, or equivalently $w(t)=t$ for all
$w\in W$. It follows that
$Z$ acts on the set $T_\ell ^{\rm reg}$ by multiplication; 
this action commutes with that of $W$ and thus defines an action of  $Z$ on
$\tl/W$. Through the bijection $P_\ell \rightarrow \tl/W$ the action of $Z$ on
$P_\ell
$ is the one deduced from its action on the extended Dynkin
diagram (see [O-W], \S 3 or [Bo], 2.3 and 4.3).

\ind Now let $\Gamma $ be a subgroup of  $Z$, and let
$G'=G/\Gamma $. We denote by $P'_\ell $ the sublattice of weights $\lambda \in
P_\ell $ such that $\lambda _{|\Gamma }=1$. The action of $\Gamma $ on $P_\ell
$ preserves $P'_\ell $; we denote by
$\Gamma\cdot\lambda
$ the orbit of a weight $\lambda $ in $P'_\ell $. 
 The Verlinde formula for $G'$ is:
$$N_\ell (G')=|\Gamma |\sum_{\lambda \in P'_\ell}|\Gamma \cdot
\lambda |^{-2g}\Bigl({|T_\ell |\over \Delta (t_\lambda )}\Bigr)^{g-1}\ .$$ 
Each term in the sum is invariant under $\Gamma $, so we may as well sum over
$P'_\ell/\Gamma  $ provided we multiply each term by $|\Gamma \cdot \lambda|
$:\global\advance\prmno by1
$$N_\ell (G')=|\Gamma
|\sum_{\lambda \in  P'_\ell/\Gamma }|\Gamma \cdot \lambda |^{1-2g}\Bigl({|T_\ell
|\over
\Delta (t_\lambda )}\Bigr)^{g-1}\ .\eqno(\the\secno.\the\prmno)$$\label{ver}
\subsec{The general case}The above formula actually applies to any semi-simple 
group
$G'=G/\Gamma $, where $G$ is a product of simply connected groups $G_1, \ldots,
G_q$ [AMW].
\ind We choose a maximal torus $T^{(i)}$ in $G_i$ for each $i$ and put
 $T=$\quad  \break$T^{(1)}\times \ldots \times $ $T^{(q)}$.
Let ${\ll }:=(\ell_1,\ldots ,\ell _q)$ be a $q$\tx uple of nonnegative integers.
We put $T_{\ll}=T^{(1)}_{\ell _1}\times \ldots \times T^{(q)}_{\ell _q}$; the
subset $T^{\rm reg}_{\ll}$ of regular elements in $T_{\ll}$ is the product of the
subsets $(T^{(i)}_{\ell _i})^{\rm reg}$.
  For each $i$, let
$P^{(i)}_{\ell _i}$ be the set of dominant weights of $T^{(i)}$ associated to
$G_i$ and $\ell _i$ as in (\ref{tl}), and let $P_{\ll}=P^{(1)}_{\ell _1}\times \ldots
\times P^{(q)}_{\ell _q}$. For $\lambda=(\lambda _1,\ldots ,\lambda _q)\in
P_{\ll}$, we put $t_\lambda =(t_{\lambda _1},\ldots ,t_{\lambda _q})\in T_{\ll}$;
this defines a bijection of
$P_{\ll}$ onto $T^{\rm reg}_{\ll}/W$.
The elements of $P_{\ll}$ are characters of  $T$, and we denote by
$P'_{\ll}$ the subset of characters which are trivial on
$\Gamma $. The group $\Gamma $ is contained in the center $Z_1\times \ldots
\times Z_q$ of $G$, which acts naturally on $P_{\ll}$ and $P'_{\ll}$.
Then\global\advance\prmno by1
$$N_{\ll} (G')=|\Gamma |\sum_{\lambda \in P'_{\ll}/\Gamma }|\Gamma \cdot
\lambda |^{1-2g}\Bigl({|T_{\ll} |\over \Delta (t_\lambda
)}\Bigr)^{g-1}\eqno(\the\secno.\the\prmno)$$\label{gen}  with $\displaystyle
\quad  {|T_{\ll}|\over \Delta (t_\lambda )}=\prod_{i=1}^q {|T_{\ell _i}|\over
\Delta_i (t_{\lambda_i})}\ $ \raise1,5pt\hbox{,}   $\displaystyle\quad  \Delta
_i(t)=\prod_{\alpha
\in R(G_i,T^{(i)})}^{}(\alpha (t)-1)$ for $t\in T^{(i)}$. 
\section{The Verlinde formula for ${\bf SO}_{\grit r}$}
\ind We now apply the previous formulas to the case $G'=SO_r$. We will
rest very much on the computations of [O-W]. We will borrow their 
notation as well as that  of [Bo].
\subsec{The case $G'=SO_{2s}$, $s\ge 3$}The root
system $R$ is of type $D_s$. 
 Let
 $(\varepsilon _1,\ldots ,\varepsilon _s)$ be the standard basis of 
${\bf R}^s$. The weight lattice $P(R)$ is spanned by the  
fundamental weights 
$$\nospacedmath\displaylines{\varpi_j=\varepsilon _1+\ldots
+\varepsilon _j\ (1\le j\le s-2),\ \cr
\varpi_{s-1}={1\over 2}(\varepsilon _1+\ldots +\varepsilon
_{s-1}-\varepsilon _s)\ ,\ \varpi_s={1\over 2}(\varepsilon _1+\ldots +\varepsilon
_{s-1}+\varepsilon _s)\ .}
$$For $\lambda \in P(R)$, we write 
$\displaystyle \lambda +\rho =\sum_i t_i\varpi_i=\sum_iu_i\varepsilon _i$
with
$$\nospacedmath\displaylines{ u_1=t_1+\ldots +t_{s-2}+{1\over
2}(t_{s-1}+t_s)\quad ,\quad
\ldots
\quad ,\quad u_{s-2}=t_{s-2}+{1\over 2}(t_{s-1}+t_s)\cr
u_{s-1}={1\over 2}(t_{s-1}+t_s)\quad ,\quad u_s={1\over
2}(-t_{s-1}+t_s)\cr
(u_i\in {1\over 2}{\bf Z}\quad ,\quad  u_i-u_{i+1}\in
{\bf Z} )\ .}$$ Put $k=\ell+2s-2$. The condition
$\lambda\in P_\ell $ becomes:
$ u_1>\ldots >u_s$,  $u_1+u_2<k$ and $u_{s-1}+u_s>0$; the condition $\lambda \in
P'_\ell $ imposes moreover 
$t_{s-1}\equiv t_s$ (mod. 2), that is, $u_i\in{\bf Z}$ for each $i$.
Thus we find a bijection between
$P'_\ell $ and the subsets
$U=\{u_1,\ldots ,u_s\}$ of ${\bf Z}$ satisfying the above 
conditions.  
\ind The group $Z$ is canonically isomorphic to $P(R)/Q(R)$
(note that $R=R^\vee$ in this case); its nonzero elements are
 the classes of $\varpi_1,\varpi_{s-1}$ and $\varpi_s$. 
The nonzero element
$\gamma $ which vanishes in $SO_{2s}$ is represented by the
only  weight in this list which comes from $SO_{2s}$, namely $\varpi_1$.
It corresponds to the automorphism of the
extended Dynkin diagram which exchanges $\alpha _0 $ with $\alpha _1$ and 
$\alpha _{s-1}$ with $\alpha _s$
 (see [Bo], Table $D_l$); it acts on 
 $P_\ell $
by $\gamma (u_1,\ldots ,u_s)= (k-u_1,u_2,\ldots ,u_{s-1},-u_s)$. Thus the
subsets $U$ as above with
$u_s\ge 0$, and moreover $u_1\le {k\over 2}$ if $u_s=0$, form a system
of representatives of
$P'_\ell /\Gamma $. The corresponding orbit has one element if
$u_1={k\over 2}$ and $u_s=0$, and 2 otherwise.

\ind For  a subset $U$ corresponding to the weight $\lambda $ we have [O-W] 
$$\Delta (t_\lambda )=\Pi _{k}(U)= \prod_{1\le i<j\le s}4\sin^2{\pi \over
k}(u_i-u_j)\  4\sin^2{\pi \over
k}(u_i+u_j)\ .$$

\ind Now we restrict ourselves to the case $\ell =2$, so that $k=r=2s$. 
Put 
$V=$ $\{s, s-{1},\ldots ,0\}$. The
subsets $U$ to consider are those of the form
 $U_j:=V\moins\{j\} $  for $0\le j\le s$. We have
\indp  $\bullet\ \Pi _r(U_j)=4r^{s-1}$  for  $ 1\le j\le s-1$ by Corollary 1.7
(ii) in [O-W];
\indp $\bullet\ \Pi_r(U_0)=\Pi_r(U_s)=r^{s-1}$ by Corollary 1.7
(iii) in [O-W].

\ind We have  $|T_2 |=4r^s$ (\ref{tl}). Multiplying the terms $U_0$ and $U_s$  by
$2^{1-2g}$ and summing, we find:
$$N_2(SO_{2s})=2[(s-1).r^{g-1}+2^{1-2g}[2.(4r)^{g-1}]=r^g\ .$$

\subsec{The case $G'= SO_{2s+1}$, $s\ge 2$}Then $R$ is of type
$B_s$. Denoting again by
 $(\varepsilon _1,\ldots ,\varepsilon _s)$  the standard basis of ${\bf
R}^s$,  the weight lattice $P(R)$ is spanned by the fundamental weights
$$\varpi_1=\varepsilon _1\ \  ,\ \  \varpi_2=\varepsilon _1+\varepsilon _2
\ \  , \ \ldots\ ,\ \   \varpi_{s-1}=\varepsilon _1+\ldots +\varepsilon
_{s-1}\ \ ,\ \
\varpi_s={1\over 2}(\varepsilon _1+\ldots +\varepsilon _s)\ .
$$For $\lambda \in P(R)$, we write
$\displaystyle \lambda +\rho =\sum_i t_i\varpi_i=\sum_iu_i\varepsilon _i$
with
$$ u_1=t_1+\ldots +t_{s-1}+{t_s\over 2}\quad ,\quad \ldots
\quad ,\quad u_{s-1}=t_{s-1}+{t_s\over 2}\quad ,\quad u_{s}={t_s\over 2}\ ,$$
with $u_i\in {1\over 2}{\bf Z}$ and $u_i-u_{i+1}\in {\bf Z}$ for each $i$.
Put $k=\ell +2s-1$. The condition $\lambda \in P_\ell $ becomes $u_1>\ldots >u_s>0$
and 
$u_1+u_2<k$. Since $\varpi_s$ is the only fundamental weight which
does not come from $SO_{2s+1}$, the condition $\lambda \in P'_\ell $ is equivalent to
$t_s$ odd, that is, $u_s\in {\bf Z}+{1\over 2}$.
Thus we find a bijection between $P'_\ell $ and the subsets
$U=\{u_1,\ldots ,u_s\}$ of ${\bf Z}+{1\over 2}$ satisfying 
$$u_1>\ldots >u_s>0\quad , \quad  u_1+u_2<k\ .$$
The non-trivial element
$\gamma $ of
$\Gamma$ acts on $P_\ell $ by $\gamma (u_1,\ldots ,u_s)=
(k-u_1,u_2,\ldots ,u_s)$. Thus the elements $U$ as above with $u_1\le {k\over
2}$ form a system of representatives of $P'_\ell /\Gamma $. The corresponding
orbit has one element if $u_1={k\over 2}$ and 2 otherwise.

\ind 
For  a subset $U$ corresponding to the weight $\lambda $ we have [O-W]
$$\Delta (t_\lambda )=\Phi _{r}(U)= \prod_{1\le i<j\le s}4\sin^2{\pi \over
r}(u_i-u_j)\  4\sin^2{\pi \over
r}(u_i+u_j)\prod_{i=1}^s4\sin^2{\pi \over
r}u_i\ .$$

\ind Now we restrict ourselves to the case $\ell =2$, so that $k=r=2s+1$. 
Put 
$V=$ $\{s+{1\over 2}, s-{1\over 2},\ldots ,{1\over 2}\}$. The
subsets $U$ to consider are 
  the subsets $U_j:=V\moins\{j+{1\over 2}\} $ for $0\le j\le s$.
We have
\ind  $\bullet\ \Phi _r(U_j)=4r^{s-1}$  for  $ 0\le j\le s-1$ by Corollary 1.9 (ii) in
[O-W];
\ind  
 $\bullet\ \Phi
_r(U_s)=r^{s-1}$ by Corollary 1.9 (ii) in [O-W].
\ind We have again $|T_2 |=4r^s$ (\ref{tl}).
 Multiplying the term $U_s$  by $2^{1-2g}$
and summing, we find:
$$N_2(SO_{2s+1})=
2[s.r^{g-1}+2^{1-2g}(4r)^{g-1}]=r^g\ .$$

\subsec{The case $G'=SO_3$}In that case $G=SL_2$ has a unique fundamental weight
$\rho $, and a unique positive root $\theta =2\rho $. The Dynkin index of the
standard representation of $SO_3$ is $4$, so we want to compute $N_4(SO_3)$.
We have $|T_4|=12$ (\ref{tl}). The set $P_4$ contains the weights
$k\rho$ with $0\le k\le 4$; the weights with
$k$ even come from $SO_3$, and $\Gamma $ exchanges $k\rho $ and
$(4-k)\rho$. Thus a system of representatives of $P'_4/\Gamma $ is 
$\{0,2\rho \}$, with $|\Gamma \cdot 0 |=2$ and $|\Gamma \cdot
2\rho  |=1$. Formula   (\ref{ver}) gives:
$$N_2(SO_3)=2.[2^{1-2g}12^{g-1}+3^{g-1}]=3^g\ . $$

\subsec{The case $G'=SO_4$}In that case $G=SL_2\times SL_2$ and the
nontrivial element of $\Gamma$ is $(-I,-I)$. The Dynkin index of the standard
representation of $SO_4$ is $(2,2)$. We have $|T_2|=8$ for $SL_2$,
hence $|T_{(2,2)}|=8^2$. The set $P_{(2,2)}$ contains the weights
$(k\rho ,l\rho) $ with $0\le k,l\le 2$, and $P'_{(2,2)}$ is defined by
the condition $k\equiv l$ (mod. 2). The element $(-I,-I)$ exchanges
$(k\rho ,l\rho )$ with $((2-k)\rho ,(2-l)\rho) $. Thus
$P'_{(2,2)}/\Gamma $ consists of the classes of $(0,0)$, $(0,2\rho )$
and
$(\rho ,\rho )$, the latter being the only one with a nontrivial stabilizer. Formula
(\ref{gen}) gives 
$$N_{(2,2)}(SO_4)=2\,[\,2\cdot 2^{1-2g}\cdot 4^{2g-2}+ 2^{2g-2}\,]=4^g\ .$$
\ind Therefore for each $r\ge 3$ we have obtained $\dim H^0(\mso,\lso)=r^g$. This
achieves the proof of Proposition \ref{main}, and therefore of the Theorem stated
in the introduction.
\section{The moduli space $\hbox{\cag\char77}_{{\bf 
Sp_{\grit 2r}}}$}
\subsection Let $r$ be an integer $\ge 1$. The  space
$\msp{r}$ is the
 moduli space  of (semi-stable) symplectic bundles, that is pairs $(E,\varphi )$
where $E$ is a semi-stable\note{By the same argument as in the orthogonal case
(footnote 1), a symplectic bundle $(E,\varphi )$ is semi-stable
if and only if $E$ is semi-stable as a vector bundle.} vector bundle of rank
$2r$ and trivial determinant  and
$\varphi :\ext^2E\rightarrow {\cal O}_C$ a non-de\-ge\-nerate alternate form.
It is connected. To alleviate the notation we will
denote it by ${\cal M}_r$. The determinant bundle  ${\cal L}_r$ generates
$\Pic({\cal M}_r)$ ([K-N], [L-S]).
\ind To describe the ``strange duality" in an intrinsic way we need a variant of
this space, namely  the  moduli space
${\cal M}'_r$ of semi-stable vector bundles $F$ of rank
$2r$ and determinant $K_C^r$, endowed with a symplectic form $\psi 
:\ext^2F\rightarrow K_C$. If $\kappa $ is a theta-characteristic on $C$, the map
$E\mapsto E\otimes \kappa $ induces an isomorphism ${\cal M}_r\iso {\cal
M}'_r$. 
We denote by  ${\cal L}'_r$ the line bundle corresponding to ${\cal L}_r$ under
any of these isomorphisms. 

\ind Similarly, we will consider for $t$ even the  moduli space ${\cal
M}'_{SO_{t}}$ of semi-stable vector bundles $E$ of rank
$t$ and determinant $K_C^{t/2}$, endowed with a quadratic form $q
:\sym^2E\rightarrow K_C$. It has two components ${\cal M}'^\pm_{SO_{t}}$
depending on the parity of $h^0(E)$; if $\kappa $ is a theta-characteristic on $C$,
the map $E\mapsto E\otimes \kappa $ induces  isomorphisms ${\cal
M}_{SO_{t}}^\pm\iso {\cal M}'^{\pm}_{SO_{t}}$ (\ref{w2}). The space ${\cal
M}'^{+}_{SO_{t}}$ carries a canonical Weil divisor, the reduced subvariety
$${\cal D} =\{(E,q)\in \ms{t}'^{+}\ |\ H^0(C,E )\not=0\}\ ;$$  $2{\cal D} $ is a
Cartier divisor, defined by a section of the generator ${\cal L}'_{SO_{t}}$ of
$\Pic(\ms{t}'^{+})$ ([L-S],
\S 7).

\subsec{The strange duality for symplectic bundles}Let $r,s$ be integers $\ge 2$,
and $t=4rs$.  Consider the map
$$\pi :{\cal M}_r\times {\cal M}'_s\longrightarrow \ms{t}'$$which
maps $((E,\varphi ),(F,\psi ) $ to $(E\otimes F,\varphi \otimes \psi)$.
Since ${\cal M}_r$ is connected and contains the trivial bundle ${\cal O}^{2r}$ with
the standard symplectic form, the image lands in $\ms{t}'^{+}$. 
\ind For $(E,\varphi)\in{\cal M}_r$, the pull back of $\ls{t}$ to
$\{(E,\varphi)\}\times {\cal M}'_s$ is the line bundle associated to $2r$ times the
standard representation, that is  ${\cal L}'^{2r}_s$; similarly its pull back to
${\cal M}_r\times
\{(F,\psi )\}$, for $(F,\psi )\in{\cal M}'_s$, is ${\cal L}_r^{2s}$. It follows that
 $$ \pi ^*\ls{t}\cong
{\cal L}_r^{2s}\boxtimes{\cal L}'^{2r}_s\ .$$
\ind If $\kappa $ is a theta-characteristic on $C$ with $h^0(\kappa )=0$,
we have 
$\pi ({\cal O}_C^{2r},\kappa ^{2s})\notin {\cal D} $ (${\cal O}_C^{2r}$ and 
$\kappa ^{2s}$ are endowed with the standard alternate forms).  Thus $\Delta :=\pi
^*{\cal D}$ is a Weil divisor on
${\cal M}_r\times {\cal M}'_s$, whose double is a Cartier divisor defined by a
section of  $({\cal L}_r^{s}\boxtimes{\cal L}'^{r}_s)^2$; but this moduli space is
locally factorial ([S], Thm. 1.2), so that
$\Delta $ is actually a Cartier divisor, defined by a section $\delta$ of
${\cal L}_r^{s}\boxtimes{\cal L}'^{r}_s$, well-defined up to a scalar. Via the
K\"unneth isomorphism we view $\delta $ as an element of 
$H^0({\cal M}_r,{\cal L}_r^{s})\otimes H^0({\cal M}'_s,{\cal L}'^r_s)$.
The {\it strange
duality} conjecture for symplectic bundles is
\th Conjecture
\enonce The section $\delta $
 induces an isomorphism
$$H^0({\cal M}_r,{\cal L}_r^{s})^*\iso H^0({\cal M}'_s,{\cal L}'^{r}_s)\ .$$
\endth\label{sdc}
\ind If the conjecture holds, the rational map $\varphi _{{\cal L}^s_r}:{\cal
M}_r\dasharrow |{\cal L}^s_r|$ is identified through $\delta ^{\sharp}$ to the map
${\cal M}_r\dasharrow |{\cal L}'^r_s|$ given by $E\mapsto \Delta _E$, where
$\Delta _E$ is the trace of $\Delta $ on $\{E\}\times {\cal M}'_s$; set-theoretically:
$$\Delta _E=\{(F,\varphi )\in {\cal M}'_s\ |\  H^0(C,E\otimes F)\not=0\}\ .$$
\ind By [O-W], we have $\dim H^0({\cal M}_r,{\cal L}_r^{s})=\dim H^0({\cal
M}_s,{\cal L}^{r}_s)$.  Therefore the
conjecture is equivalent to:
\subsection {\it The linear system $|{\cal L}'^{s}_r|$ is spanned by the
divisors $\Delta _{E} $, for $E\in {\cal M}_r$}.\label{span}
\smallskip \ind  We now specialize to the case $s=1$. The space ${\cal
M}'_1$ is the moduli space ${\cal N}$ of semi-stable rank 2 vector bundles on $C$
with determinant $K_C$; its Picard group is geenrated by the determinant bundle
${\cal L}$. 
The conjecture becomes:
\th Conjecture
\enonce The isomorphism $\delta ^{\sharp}:H^0({\cal M}_r,{\cal L}_r)^*\iso
H^0({\cal N},{\cal L}^{r})$ identifies the map $\varphi _{{\cal L}_r}:{\cal
M}_r\dasharrow |{\cal L}_r|^*$ with the rational map $E\mapsto \Delta _E$ of
${\cal M}_r$  into $|{\cal L}|$.
\endth\label{conj}
\ind By (\ref{span}) this is equivalent to saying that the linear system $|{\cal L}^r|$
on ${\cal N}$ is spanned by the divisors $\Delta_E$ for $E\in {\cal M}_r$.
\subsection Let $G$ be a semi-stable vector bundle of rank $r$ and degree $0$. 
To $G$ is associated a divisor $\Theta _G\in |{\cal L}^r|$,
supported on the set 
$$\Theta _G=\{F\in {\cal N}\ |\ H^0(C,G\otimes F)\not= 0\}$$provided this set
is $\not= {\cal N}$ [D-N]. Put
$E=G\oplus G^*$, with the standard symplectic form. We have $\Theta _G=\Theta
_{G^*}$ by Serre duality, hence
$\Delta_E={1\over 2}\Theta_E= {1\over 2}(\Theta _G+\Theta_{G^*})=$ $\Theta
_G$; thus  conjecture \ref{conj} holds if the linear system $|{\cal L}^r|$ on ${\cal
N}$ is spanned by the divisors $\Theta _G$ for $G$ semi-stable of degree
$0$. In particular, it suffices to prove that $|{\cal L}^r|$ is spanned by the divisors
$\Theta _{L_1}+\ldots +\Theta _{L_r}$, for $L_1,\ldots ,L_r\in J$. As a
consequence of [BNR], the divisors $\Theta _L$ for $L$ in $J$ span $|{\cal L}|$, so 
{\it Conjecture {\rm (\ref{conj})} holds if the multiplication map
$m_r:\sym^rH^0({\cal N},{\cal L})\rightarrow H^0({\cal N},{\cal L}^r)$ is
surjective.} \label{span2}  
\th Proposition 
\enonce Conjecture {\rm \ref{conj}} holds in the following cases:
\indp {\rm (i)} $r=2$ and $C$ has no vanishing thetanull;
\indp {\rm (ii)} $r\ge 3g-6$ and $C$ is general enough;
\indp {\rm (iii)} $g=2$, or $g=3$ and $C$ is non-hyperelliptic. 
\endth\label{cases} 
{\it Proof} : In each case  the multiplication map
$m_r:\sym^rH^0({\cal N},{\cal L})\rightarrow H^0({\cal N},{\cal L}^r)$ is
surjective. This follows from  [B1], Prop. 2.6 c) in case (i),  and
from the explicit description of ${\cal M}_{SL_2}$ in case (iii). When $C$ is generic,
 the surjectivity of $m_r$ for $r$ even $\ge 2g-4$ follows from that of $m_2$
together with [L]. We have $H^i({\cal N},{\cal L}^{j})=0$
for  $i\ge 1$ and
$j\ge -3$ by [K-N], Thm. 2.8. By [M2] this implies that the
multiplication map
$$H^0({\cal N},{\cal L})\otimes H^0({\cal N},{\cal L}^k)\longrightarrow H^0({\cal
N},{\cal L}^{k+1})$$is surjective for $k\ge \dim{\cal N}-3=3g-6$. Together with the
previous result this implies the surjectivity of $m_r$ for $r\ge 3g-6$, and therefore
by (\ref{span2}) the Proposition.\cqfd
\th Corollary
\enonce Suppose $C$ has no vanishing thetanull. There is  a canonical
isomorphism
$|{\cal L}_{2}|^*\iso |4\Theta |^+$  which identify the maps
$\varphi _{{\cal L}_{2}}:{\cal M}_{2}\dasharrow |{\cal L}_{2}|^*$
with  
$\theta :{\cal M}_{2}\dasharrow |4\Theta |^+$ such that $\theta
(E,\varphi )=\Theta _E$.
\endth
{\it Proof} : Let $i:J^{g-1}\rightarrow {\cal N}$ be the map
$L\mapsto L\oplus \iota ^*L$. The composition
$$H^0({\cal M}_2,{\cal L}_2)^*\qfl{\delta ^{\sharp}}
H^0({\cal N},{\cal L}^{2})\qfl{i^*}H^0(J^{g-1}, {\cal O}(4\Theta ))^+$$
is an isomorphism by Prop. \ref{cases} (i) and Prop. 2.6 c) of [B1]; it maps
$\varphi _{{\cal L}_2}(E,\varphi )$ to $i^*\Delta _E$. Using Serre duality again  we
find  $i^*\Delta _E={1\over 2}(\Theta _E+\Theta_{E^*})=\Theta _E$, hence the
Corollary.\cqfd 
(4.9) {\it Remarks}$.-$ 1) The corollary does not hold if $C$ has a vanishing
thetanull: the image of $\theta$ is contained in that of $i^*$, which is a proper
subspace of $|4\Theta |^+$.
\ind 2) The analogous statement for $r\ge 3$ does not
hold: the Verlinde formula implies 
  $\dim H^0({\cal N},{\cal L}^{r})>\dim H^0(J^{g-1}, {\cal O}(2r\Theta
))^+$ for $g\ge 3$, or $g=2$ and $r\ge 4$.

\vskip2cm
\centerline{ REFERENCES} \vglue15pt\baselineskip12.8pt
\def\num#1{\smallskip\item{\hbox to\parindent{\enskip [#1]\hfill}}}
\parindent=1.35cm 
\num{AMW} A. {\pc ALEXEEV}, E. {\pc MEINRENKEN}, C. {\pc WOODWARD}:
{\sl Formulas of Verlinde type for non-simply connected groups}. Preprint
{\tt math.SG/0005047}.

\num{B1} A. {\pc BEAUVILLE}: {\sl Fibr\'es de rang $2$ sur les courbes, fibr\'e 
d\'eterminant et fonctions th\^eta II}.  Bull. Soc. math. France {\bf 119}
(1991), 259--291.
\num{B2} A. {\pc BEAUVILLE}: {\sl Conformal blocks, Fusion rings and the Verlinde
formula.} Proc. of the Hirzebruch 65 Conf. on Algebraic Geometry, Israel Math. Conf.
Proc. {\bf 9}  (1996), 75--96.
\num{B3} A. {\pc BEAUVILLE}: {\sl Vector bundles on curves and  theta functions}.
Preprint {\tt math.AG/0502179}.
\num{BLS} A. {\pc BEAUVILLE}, Y. {\pc LASZLO}, C. {\pc SORGER}: {\sl The Picard
group of the moduli of $G$\tx bundles on a
 curve}. Compositio math. {\bf 112} (1998), \no 2, 183--216.
\num{BNR} A. {\pc BEAUVILLE}, M.S. {\pc NARASIMHAN}, S. {\pc
RAMANAN}: {\sl Spectral curves and the generalised theta
divisor}. J. Reine Angew. Math. {\bf 398} (1989), 169--179. 

\num{Bo} N. {\pc BOURBAKI}: {\sl Groupes et alg\`ebres de Lie}, Chap.\ VI.
Hermann, Paris (1968).

\num{D-N} J.-M. {\pc DREZET}, M.S. {\pc NARASIMHAN}: {\sl Groupe de Picard des
vari\'et\'es de modules de fibr\'es semi-stables sur les courbes alg\'ebriques}. 
Invent. Math. {\bf 97} (1989), \no 1, 53--94.

\num{D} E. {\pc DYNKIN}: {\sl Semisimple subalgebras of semisimple Lie algebras}.
Amer. Math. Soc. Translations (II) {\bf 6} (1957), 111--244. 

\num{K-N} S. {\pc KUMAR}, M.S. {\pc NARASIMHAN}: {\sl Picard group of the
moduli spaces of $G$-bundles}. Math. Ann. {\bf 308} (1997), \no 1, 155--173. 

\num{L} Y. {\pc LASZLO}: {\sl  \`A propos de l'espace des modules de fibr\'es de
rang $2$ sur une courbe}. Math. Ann. {\bf 299} (1994), \no 4, 597--608.

\num{L-S} Y. {\pc LASZLO}, C. {\pc SORGER}:
{\sl The line bundles on the moduli of parabolic $G$-bundles over curves and their
sections}. Ann. Sci. \'Ecole Norm. Sup. (4) {\bf 30} (1997), \no 4, 499--525.

\num{M1} D. {\pc MUMFORD}: {\sl On the equations defining abelian varieties} I.
Invent. Math. {\bf 1} (1966), 287--354. 

\num{M2} D. {\pc MUMFORD}: {\sl Varieties defined by quadratic equations}. 
Questions on Algebraic Varieties (1970), 29--100;
Cremonese, Rome.

\num{O-W} W. {\pc OXBURY}, S. {\pc WILSON}: {\sl  Reciprocity laws in the Verlinde
formulae for the classical groups}. Trans. Amer. Math. Soc. {\bf 348} (1996), \no 7,
2689--2710.

\num{R} S. {\pc RAMANAN}: {\sl Orthogonal and spin bundles over
hyperelliptic curves}. Proc. Indian Acad. Sci. Math. Sci. {\bf 90} (1981), \no 2,
151--166.
\num{Se} J.-P. {\pc SERRE}: {\sl Rev\^etements \`a ramification impaire et
th\^eta-caract\'eristiques}. C. R. Acad. Sci. Paris S\'er. I Math. {\bf 311} (1990), \no
9, 547--552.
\num{S} C. {\pc SORGER}: {\sl On moduli of $G$-bundles of a curve for
exceptional $G$}. Ann. Sci. \'Ecole Norm. Sup. (4) {\bf 32} (1999), \no 1, 127--133. 
\vskip1cm

\def\pc#1{\eightrm#1\sixrm}
\hfill\vtop{\eightrm\hbox to 5cm{\hfill Arnaud {\pc BEAUVILLE}\hfill}
 \hbox to 5cm{\hfill Institut Universitaire de France\hfill}\vskip-2pt
\hbox to 5cm{\hfill \&\hfill}\vskip-2pt
 \hbox to 5cm{\hfill Laboratoire J.-A. Dieudonn\'e\hfill}
 \hbox to 5cm{\sixrm\hfill UMR 6621 du CNRS\hfill}
\hbox to 5cm{\hfill {\pc UNIVERSIT\'E DE}  {\pc NICE}\hfill}
\hbox to 5cm{\hfill  Parc Valrose\hfill}
\hbox to 5cm{\hfill F-06108 {\pc NICE} Cedex 02\hfill}}

\end